\documentclass[english]{cccconf}
\usepackage[comma,numbers,square,sort&compress]{natbib}
\usepackage{epstopdf}

\usepackage{amsmath,amssymb,amsfonts}
\usepackage{algorithm,algorithmic}
\usepackage{graphicx}
\usepackage{textcomp}
\usepackage{url}
\usepackage{color}

\usepackage{lipsum}

\usepackage{amsopn}
\DeclareMathOperator{\diag}{diag}

\usepackage{enumerate}
\usepackage{cancel}
\usepackage{mathrsfs}
\usepackage{mathdots}
\usepackage{euscript}
\usepackage{amscd}
\usepackage{placeins}

\usepackage{subcaption}

\usepackage{tikz}
\usetikzlibrary{arrows,shapes,chains}

\newtheorem{theorem}{Theorem}

\newtheorem{example}{Example}

\newtheorem{remark}{Remark}

\newcommand{\bmat}{\left[ \begin{matrix}}
	\newcommand{\emat}{\end{matrix} \right]}

\DeclareMathOperator{\trace}{tr}

\DeclareMathOperator{\E}{{\mathbb E}}
\newcommand{\Rbb}{\mathbb R}

\newcommand{\Cbb}{\mathbb C}

\newcommand{\xb}{\mathbf  x}

\newcommand{\sbf}{\mathbf  s}  

\newcommand{\wb}{\mathbf  w}

\newcommand{\bb}{\mathbf  b}
\newcommand{\db}{\mathbf  d}

\newcommand{\ub}{\mathbf  u}

\newcommand{\Gb}{\mathbf G}

\newcommand{\rhob}{\boldsymbol{\rho}}
\newcommand{\thetab}{\boldsymbol{\theta}}

\DeclareMathOperator{\range}{Range}
\DeclareMathOperator{\rank}{rank}

\newcommand{\Gscr}{\mathscr{G}}

\newcommand{\Acal}{\mathcal{A}}

\newcommand{\Ical}{\mathcal{I}}
\newcommand{\Dcal}{\mathcal{D}}
\newcommand{\Hcal}{\mathcal{H}}

\renewcommand{\d}{\mathrm{d}}

\newcommand{\SNR}{\mathrm{SNR}}

\newcommand{\F}{\mathrm{F}}

\newcommand{\srm}{\mathrm{s}}

\begin{document}

\title{Line Spectral Estimation Using a G-Filter: Atomic Norm Minimization with Multiple Output Vectors}

\author{Jiale Tang\aref{sysu},
        Bin Zhu\aref{sysu}}

\affiliation[sysu]{School of Intelligent Systems Engineering,
        Sun Yat-sen University, Shenzhen 518107, P.~R.~China
        \email{tangjle@mail2.sysu.edu.cn, zhub26@mail.sysu.edu.cn}}

\maketitle

\begin{abstract}
We propose an atomic norm minimization (ANM) estimator of frequencies in a noisy complex sinusoidal signal that integrates Georgiou’s filter bank (G-filter) with multiple output vectors (MOV). Unlike our previous work on the G-filter version of ANM which is restricted to a single filtered output vector, the proposed method in this paper uses MOV to improve data utilization and robustness of the estimate. 
The ANM problem with MOV can be reformulated as a semidefinite program thanks to a Carath\'{e}odory--Fej\'{e}r-type  decomposition for output covariance matrices of the G-filter. Numerical simulations demonstrate that the proposed approach significantly outperforms the standard ANM and the G-filter version of ANM with a single output vector in recovering the correct number of frequency components when the frequencies fall within the band(s) selected by the G-filter, particularly in the low SNR regime.
\end{abstract}


\footnotetext{This work was supported in part by Shenzhen Science and Technology
Program (Grant No. 202206193000001-20220817184157001), and the
``Hundred-Talent Program'' of Sun Yat-sen University.}

\section{Introduction}

Estimating the power distribution of a signal across the frequency domain from measurements at a finite number of time instances remains a fundamental problem in the fields of systems and control and signal processing \cite{LP15,stoica2005spectral}. In this paper, we consider recovery of frequency components from signals that can be expressed as a superposition of complex sinusoids corrupted by noise, a problem commonly referred to as \emph{line spectral estimation} or \emph{frequency estimation}, or \emph{direction-of-arrival estimation} in array processing.
The problem has a wide range of applications, including remote sensing \cite{feng2018airborne}, channel estimation in wireless communications \cite{bajwa2010compressed}, imaging \cite{borcea2002imaging}, and others \cite{Picci-Zhu-2021}.

Traditional methods in line spectral estimation face limitations: FFT suffers from ``basis mismatch'' due to the gridding operation (discretization) of the frequency domain, while subspace methods may have poor performance in the case of few measurements. Motivated by compressed sensing (see e.g., \cite{donoho2006compressed}), the \emph{atomic norm minimization} (ANM) method provides a gridless alternative which also exhibits a superior small-sample performance (see e.g., \cite{tang2013compressed,bhaskar2013atomic,li2015off,ZHU2026112759}). In particular, the standard ANM exploits the classic Carath\'{e}odory--Fej\'{e}r (C--F) decomposition to represent spectral lines through a positive semidefinite Toeplitz matrix, which yields a tractable (convex) semidefinite program.

However, the standard ANM treats all frequencies uniformly and lacks frequency selectivity, making it unable to exploit \emph{a priori} knowledge about underlying frequencies when such knowledge is available in certain applications. To overcome this, we employ a class of stable linear filters advocated by Georgiou and coauthors which we refer to as ``G-filter'' (see e.g., \cite{georgiou2000signal,amini2006tunable}) to select the frequency band of interest. The use of a G-filter is also justified by a generalized C--F-type decomposition which holds for a broader class of output covariance matrices than Toeplitz matrices that are compatible with the filter structure.

While the G-filter has been used primarily for estimating rational spectral densities (see e.g.,~\cite{BGL-THREE-00,Georgiou-L-03,FRT-11,Z-14rat,cao2024spectral}), our recent works \cite{ZHU2026112759,Zhu_CDC2025} have adopted a different viewpoint by leveraging the framework to embed frequency selectivity into ANM for line spectral estimation.

Despite the fact that our previous method is effective with only one output vector of the G-filter, its performance worsens in the regime of low signal-to-noise ratios (SNRs). 
To address this, in this paper we propose an improved approach that utilizes multiple output vectors (MOV) of the G-filter. We reformulate the ANM problem within this MOV setting, ensuring stable recovery of the frequencies in a wide range of SNRs. Through the use of MOV, the proposed method enhances robustness of estimation in challenging cases, particularly when the SNR is low or two frequencies are close to each other. In addition, we also illustrate a design of the G-filter for frequency estimation across multiple passbands, thereby broadening the applicability of our approach.

The rest of this paper is organized as follows. A general signal model using a G-filter with MOV is described in Section~\ref{sec_Signal_model}. 
The C--F-type decomposition for output covariance matrices of a G-filter is reviewed in Section~\ref{sec_ANM} and an ANM problem for frequency estimation in the MOV case is discussed. Detailed numerical simulations are presented in Section~\ref{sec_sims}. Finally, Section~\ref{sec_conc} draws conclusions.

\section{Signal model and problem statement}\label{sec_Signal_model}

Consider a discrete signal $y(t)$ modeled as a superposition of $m$ cisoids (complex sinusoids) corrupted by additive complex white noise $w(t)$ and measured at $L$ time instances:
\begin{equation}\label{signal_model}
y(t) = s(t) + w(t) = \sum_{k=1}^m a_k \, e^{i \theta_k t} + w(t), \quad t = 0, \dots, L-1,
\end{equation}
where $a_k \in \mathbb{C}$ and $\theta_k \in [0, 2\pi)$ denote an amplitude and an unknown angular frequency, respectively. The number of cisoidal components $m$ is unknown.

The G-filter, originally introduced in \cite{georgiou2000signal}, is defined by a state  equation
\begin{equation}\label{filter_bank_rewrite}  
\xb(t) = A\xb(t-1) + \bb y(t), \quad t \in \mathbb{Z},  
\end{equation}  
where the input $y(t)$ of the filter is in \eqref{signal_model} and the output $\xb(t)$ is $\Cbb^n$-valued. The state matrix $A \in \mathbb{C}^{n \times n}$ is discrete-time (Schur) stable, which means that all its eigenvalues have modulus strictly less than one. The vector $\bb \in \mathbb{C}^n$ forms a \emph{reachable} pair with $A$, that is, the condition $\rank \begin{bmatrix} \bb & A\bb & \cdots & A^{n-1}\bb \end{bmatrix} = n$ holds.  
The transfer function of the system/filter \eqref{filter_bank_rewrite} can be expressed as
\begin{equation}\label{trans_func_filter_bank_rewrite}  
G(z) = \begin{bmatrix} g_1(z) \\ \vdots \\ g_n(z) \end{bmatrix} = z(zI - A)^{-1} \bb, 
\end{equation}  
where $z$ can be interpreted as the forward shift operator $\xb(t) \mapsto \xb(t+1)$. 

\begin{example}[A bandpass filter bank]\label{trans_func_filter_bank}
	Take $A$ to be the $n\times n$ Jordan block $J_1$ with a complex constant $p_1$ on the main diagonal such that $|p_1| < 1$, that is,
	\begin{equation}\label{delay_filter}
	J_1=\bmat p_1 & 1 & \cdots & 0 \\
	\vdots & \ddots & \ddots & \vdots \\
	0 & \cdots & p_1 & 1 \\
	0 & \cdots & 0 & p_1\emat
	\quad \text{and} \quad
	\bb =\bmat 0 \\ \vdots \\ 0 \\ 1 \emat.
	\end{equation}
Then the $k$-th component of the transfer function is
\begin{equation}
	g_k(z) = \frac{z}{(z - p_1)^{n - k + 1}}, \quad k=1, 2, \dots, n.
\end{equation}
An interesting special case is obtained by setting the parameter $p_1 = 0$, in which case
\begin{equation}
	g_k(z) = z^{-n+k}, \quad k=1, 2, \dots, n,
\end{equation}
which represents a delay filter bank, the standard model for frequency estimation in the signal-processing literature. When applied to the input signal \eqref{signal_model}, the delay filter bank stacks the scalar measurements $\{y(t)\}_{t=0}^{L-1}$ into vectors of size $n$ without further processing.
\end{example}

\begin{figure}
	\centering
	\includegraphics[width=0.35\textwidth]{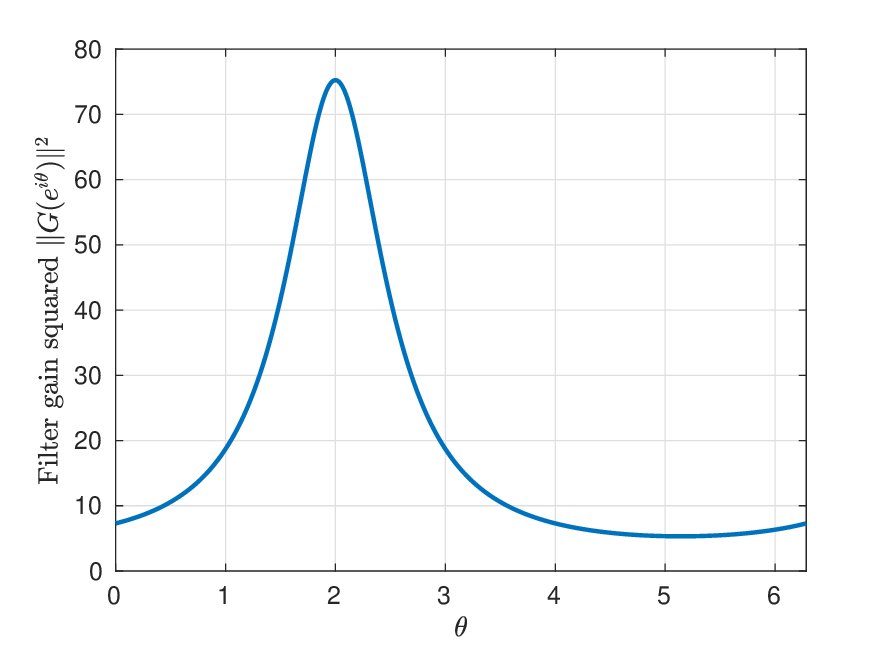}
	\caption{The squared gain $\|G^{(1)}(e^{i\theta})\|^2$ of a G-filter of size $n=20$ versus the frequency $\theta\in [0, 2\pi)$. The filter parameters $(A^{(1)}, \bb^{(1)})$ are constructed as per Example~\ref{trans_func_filter_bank} with a repeated pole $p_1=0.58e^{i2}$.}
	\label{fig:filter-gain-one}
\end{figure}

Example~\ref{trans_func_filter_bank} shows that via G-filtering the signal model for frequency estimation becomes formally more general. In addition to mere generality, the G-filter can be used to perform \emph{band selection}, i.e., to amplify the frequency content of the input $y(t)$ in certain subintervals of the frequency domain $[0, 2\pi)$, which can improve the accuracy of frequency estimates when the true frequencies reside in that band.
Note that the passbands represent our \emph{a priori} knowledge about the frequencies in the cisoids, see e.g., \cite{ZHU2026112759,amini2006tunable,georgiou2000signal,yang2018frequency}.
For instance, the G-filter in Example~\ref{trans_func_filter_bank} with a repeated pole at $p_1=\rho e^{i\varphi}$ of multiplicity $n$ can be used to select a frequency band $[\theta_{\ell}, \theta_{u}]$ near $\varphi$ that we set. Following a design procedure in \cite[p.~2667]{amini2006tunable}, we take $p_1 = 0.58e^{i2}$ and the filter size $n=20$, and the resulting G-filter, denoted by $G^{(1)}$, selects the frequency band $\mathcal{I}_1 = [1.75, 2.25]$. Indeed, such a band selection can be seen from Fig.~\ref{fig:filter-gain-one} in which we plot the squared norm $\|G^{(1)}(e^{i\theta})\|^2$ which represents a gain of the G-filter. Clearly, the gain is large when $\theta$ is near $2$ which corresponds to the designed pole $p_1$.

\begin{figure}
	\centering
	\includegraphics[width=0.35\textwidth]{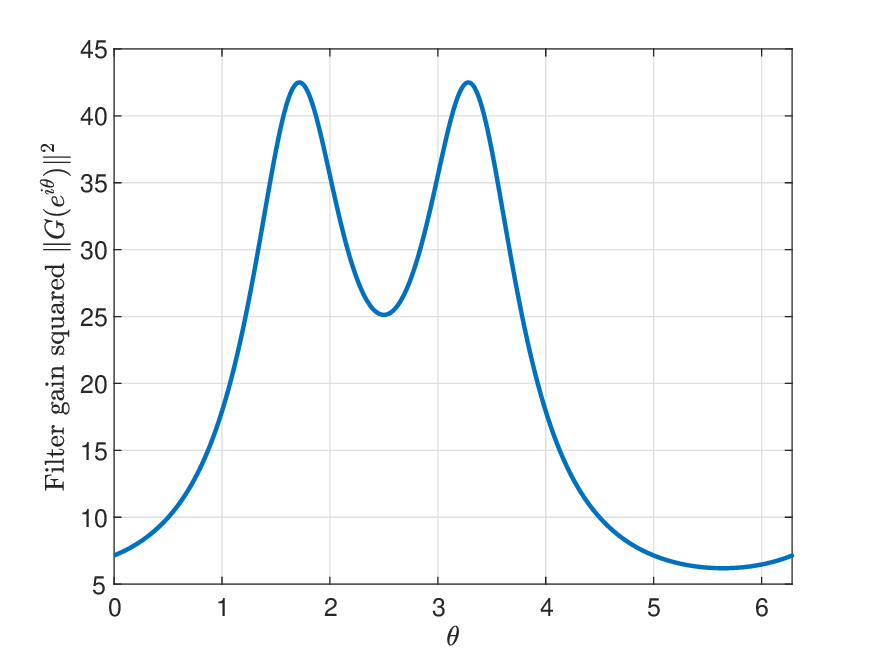}
	\caption{The squared gain $\|G^{(2)}(e^{i\theta})\|^2$ of a G-filter in Example~\ref{filter-gaintwo-example} of size $n=20$ versus the frequency $\theta\in [0, 2\pi)$. The filter parameters $(A^{(2)}, \bb^{(2)})$ are constructed with two poles $p_{k}=0.58e^{i\varphi_k}$, $k = 1,2$ where $\varphi_1 = 1.7$ and $\varphi_2 = 3.3$, respectively, {each having multiplicity $n/2=10$}.}
	\label{fig:filter-gain-two}
\end{figure}

\begin{example}[G-filter with multiple passbands]\label{filter-gaintwo-example}
The G-filter constructed in Example~\ref{trans_func_filter_bank} can be expanded to accommodate different poles so that multiple passbands can be selected. We illustrate the case of two poles $p_1=\rho_1 e^{i\varphi_1}$ and $p_2=\rho_2 e^{i\varphi_2}$ here and make equal the sizes of their corresponding Jordan blocks $J_1$ and $J_2$. The filter parameters $(A,\bb)$ are taken as a block diagonal $A = \diag(J_1,J_2)$ and $\bb = \bmat \bb_1;\bb_2 \emat$ where each $\bb_k$ has a compatible size with $J_k$ for $k=1, 2$ and the semicolon denotes vertical stacking. To illustrate the design, we choose $p_1 = 0.58e^{i1.7}$ and $p_2 = 0.58e^{i3.3}$, which results in a G-filter, denoted by $G^{(2)}$, that selects two frequency bands $\mathcal{I}_1 = [1.45, 1.95]$ and $\mathcal{I}_2 = [3.05, 3.55]$. This is illustrated in Fig.~\ref{fig:filter-gain-two} which shows that the squared gain of the G-filter exhibits two distinct peaks centered approximately at $\theta = 1.7$ and $\theta = 3.3$. 
\end{example}

\begin{remark}
Suppose that we are given the filter parameters $(\tilde A, \tilde \bb)$ such as those in Examples~\ref{trans_func_filter_bank} and \ref{filter-gaintwo-example}.
For technical reasons, we must enforce the normalization condition $A A^* + \bb \bb^* = I$, see \cite{amini2006tunable} and also \cite[Sec.~9.1]{ZHU2026112759}. This is achieved by solving the associated Lyapunov equation $E - \tilde A E \tilde{A}^* = \tilde\bb \tilde\bb^*$ for the unique $E$, followed by a similarity transform $(A, \bb) = (E^{-1/2} \tilde A E^{1/2}, E^{-1/2}\tilde\bb)$ where $E^{1/2}$ is an arbitrary matricial square root of $E$.
\end{remark}

With reference to \cite[Sec.~6]{ZHU2026112759}, the input-output relationship of the G-filter \eqref{filter_bank_rewrite} can be rewritten as
\begin{equation}\label{state_plus_noise}
	\begin{aligned}
		\xb(t) & = G(z) y(t) := \tilde{\sbf}(t) + \tilde{\wb}(t).
	\end{aligned}
\end{equation}
In the above expression, the signal part is given by
\begin{equation}\label{state_atomic_decomp}
	\tilde{\sbf}(t):=\sum_{k=1}^m G(e^{i\theta_k}) c_k(t)=\Gb(\thetab) \mathbf{c}(t),
\end{equation}
where $c_k(t) = a_k e^{i\theta_k t}$ is a filtered amplitude, $\mathbf{c}(t) = \bmat c_1(t) & \cdots & c_m(t) \emat^\top \in\Cbb^m$
and $\Gb(\thetab) := \bmat G(e^{i\theta_1}) & \cdots & G(e^{i\theta_m}) \emat\in\Cbb^{n\times m}$. The filtered noise process $\tilde\wb(t)$ is $\Cbb^n$-valued. Assume that we can measure the output $\xb(t)$ of the G-filter at $L_{\xb}$ successive  time instances, collected in a matrix $X = \bmat \xb(0) & \cdots & \xb(L_{\xb}-1)\emat \in\Cbb^{n \times L_{\xb}}$. Then the matricial form of \eqref{state_plus_noise} is just
\begin{equation}\label{model_mat}
    \begin{split}
	X & = \tilde S + \tilde W = \Gb(\thetab) C + \tilde W,
    \end{split}
\end{equation} 
where $\tilde S := \bmat \tilde{\sbf}(0) & \cdots & \tilde{\sbf}(L_{\xb}-1) \emat \in\Cbb^{n \times L_{\xb}}$, the noise matrix $\tilde W := \bmat \tilde\wb(0) & \cdots &\tilde\wb(L_{\xb}-1) \emat \in\Cbb^{n \times L_{\xb}}$, and $C := \bmat \mathbf{c}(0) & \cdots  & \mathbf{c}(L_{\xb}-1) \emat \in\Cbb^{m \times L_{\xb}}$. 

Our core task is to estimate the integer $m$ and the unknown frequencies $\{\theta_k\}_{k=1}^m$ from the filtered output matrix $X$. In the remaining part of this paper, we assume that the number of cisoids $m$ is strictly less than the filter size $n$ (i.e., $m<n$) also for technical reasons. 

\begin{remark}[Filtering finite-length signals]\label{rem:input_filtering_revised}
	The cisoidal signal $y(t)$ of length $L$ in \eqref{signal_model} must go through the G-filter \eqref{filter_bank_rewrite} in order to produce the output matrix $X$.
    Assume a zero initial condition for the state vector ($\xb(-1)=\mathbf{0}$). Following \cite{amini2006tunable}, we mitigate transient effects of the initial condition by discarding the first $L_{\srm}$ filtered samples, where $L_{\srm}$ is selected to satisfy $\|A^{L_{\srm}}\| < \varepsilon$ for a given tolerance $\varepsilon > 0$. In this way, the number of output vectors that remain is $L_{\xb} = L-L_{\srm}$, and the time instances in the matricial model \eqref{model_mat} are obtained by a shift of the origin along the time axis.
    We mention that the case of a single output vector ($L_{\xb}=1$) has been treated in \cite{ZHU2026112759}.
\end{remark}


\section{Atomic norm minimization approach}\label{sec_ANM}
In this section, we sketch the atomic norm minimization framework designed for multiple output vectors of the G-filter. We begin by defining the atomic set and the corresponding atomic norm, followed by a review of the generalized Carath\'{e}odory--Fej\'{e}r-type decomposition. Finally, we transform the atomic norm minimization problem into a tractable semidefinite program.

\subsection{Atomic norm for MOV of the G-filter}
We now define the atomic norm of the signal part of the output  matrix $\tilde S = \Gb(\thetab)C$ in \eqref{model_mat}.
First we normalize the rows of the matrix $C$ and rewrite it as
\begin{equation}
    C=\bmat \text{---} & \tilde{c}_1 \mathbf{d}^*_1 & \text{---} \\ 
    & \vdots & \\ 
    \text{---} & \tilde{c}_m \mathbf{d}^*_m & \text{---}\emat   
\end{equation}
where for each $k=1,\dots, m$, the column vector $\mathbf{d}_k\in\Cbb^{L_{\xb}}$ is such that $\|\mathbf{d}_k\| = 1$ and $\tilde c_k$ is the length of the $k$-th row of $C$. Then the matrix $\tilde S$ can be expressed as 
\begin{equation}\label{atomic_decomp}
	\tilde S =  \sum_{k=1}^m \tilde{c}_k G(e^{i\theta_k})\mathbf{d}_k^*.
\end{equation}
In this way, the noiseless output matrix $\tilde{S}$ of the G-filter is equal to a linear combination of elements from the atomic set
\begin{equation}\label{dictionary_MMV}
	\Acal := \{\Gscr(e^{i\theta})= G(e^{i\theta})\mathbf{d}^* : \theta\in [0, 2\pi),\ \|\mathbf{d}\| = 1\}.
\end{equation} 
Every member of the set $\Acal$ is called an ``atom''. The linear combination of rank-one matrices in \eqref{atomic_decomp} is called an \emph{atomic decomposition} where each unknown frequency $\theta_k$ is encoded in the selected atom $\Gscr(e^{i\theta})$. Inspired by the literature on compressed sensing, the atomic norm of $\tilde S$ is $\|\tilde S\|_{\Acal} :=$
\begin{equation}\label{atomic_norm_mmv}
	\begin{aligned}
		 & \inf_{\substack{\tilde{c}_k\in\Rbb,\, \theta_k,\\ \db_k\in\Cbb^{L_{\xb}}}} \left\{ \sum_{k} |\tilde{c}_k| \|G(e^{i\theta_k})\| : \tilde S = \sum_{k} \tilde{c}_k G(e^{i\theta_k}) \db_k^*, \right. \\
		& \left. \qquad \qquad \theta_k\in [0, 2\pi),\ \tilde{c}_k\neq 0,\ \|\db_k\|=1 \right\}.
	\end{aligned}
\end{equation}
This formulation can be interpreted as a weighted $\ell_1$ norm of the signal matrix $\tilde S$ that encourages as few terms as possible in the atomic decomposition \eqref{atomic_decomp}. A distinct feature here lies in the construction of the atomic set \eqref{dictionary_MMV} which depends continuously on the parameters $\theta$ and $\db$. Such a feature can increase the accuracy and flexibility of frequency estimation by avoiding the issue of ``grid mismatch'' inherent in traditional discretization of the frequency domain $[0, 2\pi)$, e.g., in the FFT method.

\subsection{C--F-type decomposition of output covariance matrices}\label{subsec_C-F dep}
Let $\E$ be the mathematical expectation operator. The output covariance matrix of the G-filter $\Sigma :=\E \bmat \xb(t)\xb(t)^*\emat$ admits the following integral representation
\begin{equation}\label{state_cov_mat}
	\Sigma  = \frac{1}{2\pi} \int_{0}^{2\pi} G(e^{i\theta})\, \d\mu_y(\theta)\, G^*(e^{i\theta}),
\end{equation}
where $\d\mu_y$ is the (nonnegative) spectral measure of the input signal $y$. If $\d\mu_y$ is \emph{absolutely continuous} \cite{rudin1987real} with respect to the Lebesgue measure, then it can be written as $\d\mu_y(\theta) = \Phi_y(\theta) \d\theta$ where $\Phi_y$ is the power spectral density. The output covariance matrix $\Sigma$ is by definition positive semidefinite. In addition, it resides in a vector space dictated by the G-filter which is explained next. 
We denote by $\Dcal$ and $\Hcal_n$ the real linear spaces consisting of \emph{signed} measures on $[0, 2\pi)$ and Hermitian matrices of size $n$, respectively. Then let us define the linear operator
\begin{equation}
	\begin{aligned}
		\Gamma : \Dcal & \to \Hcal_n \\
		\d\mu & \mapsto \frac{1}{2\pi} \int_{0}^{2\pi} G(e^{i\theta})\, \d\mu(\theta)\, G^*(e^{i\theta}).
	\end{aligned}
\end{equation}
Equation \eqref{state_cov_mat} implies that $\Sigma\in \range\Gamma$, and the latter is a linear subspace of $\Hcal_n$. An equivalent characterization, as shown in \cite[Prop.~3.2]{ferrante2012maximum}, is given by the equality
\begin{equation}\label{equal_constraint}
	(I-\Pi_{\bb})(\Sigma-A\Sigma A^*)(I-\Pi_{\bb}) = O,
\end{equation}
where $(A,\bb)$ are the parameters of the G-filter in \eqref{filter_bank_rewrite}, and $\Pi_{\bb}:=\bb\bb^*/(\bb^*\bb)$ is the rank-one projection matrix associated with the vector $\bb$.

\begin{theorem}[C--F-type decomposition \cite{georgiou2000signal}] \label{thm:C-F_type_decomp}
Let $\Sigma$ satisfy \eqref{state_cov_mat} for some nonnegative measure $\d\mu_y$ and have rank $r<n$. Then  it admits a unique decomposition
\begin{equation}\label{C-F-type-decomp}
    \Sigma = 
    \Gb(\thetab) \diag({\boldsymbol{\rho}}) \Gb(\thetab)^*, 
\end{equation}
where, $\rhob = \bmat \rho_1 &  \cdots & \rho_r \emat^\top \in\Rbb^r$ has each component $\rho_k > 0$, $\diag({\boldsymbol{\rho}})$ is a diagonal matrix with diagonal entries in $\rhob$, and $\thetab = \bmat \theta_1 &  \cdots &\theta_r \emat^\top$ consists of distinct frequencies $\theta_k \in [0, 2\pi)$ for $k=1,\dots,r$.

\end{theorem}

As detailed in \cite[Prop.~2]{georgiou2000signal} (see also \cite[Example~3]{ZHU2026112759}), this C--F-type decomposition can be computed numerically. More precisely, let $U_{r+1:n} :=\bmat \ub_{r+1} & \cdots & \ub_n\emat$ be  the matrix of (normalized) eigenvectors of the rank-deficient $\Sigma$ corresponding to the zero eigenvalue. Then the frequencies $\{\theta_k\}_{k=1}^r$ are the distinct points of minimum
of the nonnegative function
\begin{equation}\label{rat_func_sym}
    d(e^{i\theta}, e^{-i\theta}) = G^*(e^{i\theta}) U_{r+1:n}\, U_{r+1:n}^* G(e^{i\theta})
\end{equation}
on $[0, 2\pi)$.

\subsection{Reformulation via semidefinite programming}

Direct calculation of the atomic norm from the definition in \eqref{atomic_norm_mmv} seems rather nontrivial. 
Fortunately, following previous works \cite{tang2013compressed,li2015off, ZHU2026112759}, we find that the evaluation of $\|\tilde S\|_{\Acal}$ can be reformulated as a \emph{semidefinite program} (SDP), and the latter can be efficiently solved using established convex optimization techniques \cite{bv_cvxbook}. The formal result is stated next.

\begin{theorem}\label{thm_AN_noiseless}
        Given an output matrix $\tilde S\in\Cbb^{n\times L_{\xb}}$ of the G-filter when the input is the noiseless cisoidal signal $s(t)$ in \eqref{signal_model}, let $p$ be the optimal value of the semidefinite program    
	\begin{subequations}\label{AN_semidef_program}
		\begin{align}
		& \underset{\substack{Z \in \Hcal_{L_{\xb}},\\ \Sigma\in\Hcal_n}}{\text{minimize}}
		& & \frac{1}{2} (\trace Z + \trace\Sigma) \label{obj_noiseless} \\
		& \text{subject to}
		& & \bmat Z &\tilde S^* \\ \tilde S&\Sigma \emat \geq 0 \label{LMI_constraint} \\
			& & &  \quad \text{and}\quad \eqref{equal_constraint},
		\end{align}
	\end{subequations}
    where $\Hcal_{L_{\xb}}$ denotes the set of Hermitian matrices of size $L_{
    \xb}$. 
    Then the atomic norm $\|\tilde S\|_\Acal = p$.
\end{theorem}

The proof of the above theorem can be constructed along the lines of \cite[Theorem~2]{ZHU2026112759} and is omitted here. The solution of the SDP \eqref{AN_semidef_program} naturally leads to estimates of the frequency components in the input signal $s(t)$. In fact, from the optimal $\hat{\Sigma}$, we can extract the frequency estimates
 $\hat{\thetab}=(\hat\theta_1,\dots,\hat\theta_{\hat{r}})$ by applying the C--F-type decomposition in Theorem~\ref{thm:C-F_type_decomp} where $\hat r$ is the \emph{numerical rank} of $\hat{\Sigma}$.
 
The general case of frequency estimation in the presence of noise can be handled similarly. Again, we assume that we have access to the noisy output matrix $X\in\Cbb^{n\times L_{\xb}}$. Using $\|\tilde S\|_{\Acal}$ as a \emph{regularization} term, we formulate the denoising problem as
\begin{equation}\label{ANM_noisy}
	\underset{\tilde S \in \Cbb^{n \times L_{\xb}}}{\text{minimize}} \quad
	 \frac{1}{2} \|X- \tilde S \|_\F^2 + 
	\lambda \|\tilde S\|_{\Acal},
\end{equation} 
where $\lambda>0$ is a regularization parameter and $\|\cdot\|_\F$ denotes the Frobenius norm of a matrix. In view of Theorem~\ref{thm_AN_noiseless}, such a problem also admits an SDP description:
\begin{subequations}\label{noisy_SDP}
	\begin{align}
	& \qquad \underset{\substack{Z \in \Hcal_{L_{\xb}},\, \Sigma \in\Hcal_n \\ \tilde S \in \Cbb^{n \times L_{\xb}}}}{\text{minimize}}
	& & \frac{1}{2} \|X- \tilde S\|_\F^2 + \frac{\lambda}{2} (\trace Z + \trace\Sigma) \label{obj_noisy} \\
	& \qquad\quad  \text{subject to}
	& & \eqref{LMI_constraint} \quad \text{and}\quad \eqref{equal_constraint}.
	\end{align}
\end{subequations}
Once the SDP \eqref{noisy_SDP} is solved, the frequency estimates can be obtained from the optimal $\hat\Sigma$ as described above in the noiseless case.

\section{Simulations}\label{sec_sims}

In this section, we conduct numerical simulations to evaluate our frequency estimation method, abbreviated as ``G-filter MANM'', against the G-filter version of ANM in the single output vector (i.e., $L_{\xb}=1$) case \cite{ZHU2026112759}, abbreviated as ``G-filter SANM'', and the standard ANM \cite{bhaskar2013atomic}. Specific implementation details are described below.

\subsection{Choice of the regularization parameter $\lambda$}\label{subsec_choice_lambda}
According to \cite[Theorem~3]{li2015off}, to ensure a stable recovery of the signal component $\tilde S$ one must choose $\lambda\geq \E\|\tilde W\|_{\Acal}^*$ where $\|\cdot\|_{\Acal}^*$ denotes the \emph{dual} atomic norm. Although $\E\|\tilde W\|_{\Acal}^*$ can be estimated in the special case of a delay filter bank, extending that technique to the more general atomic set \eqref{dictionary_MMV} seems nontrivial and is left as future work. 

In the simulations presented in the next subsection, we use the heuristic value of $\lambda$ for a delay filter bank \cite{li2015off}. First, we compute an estimate $\hat\sigma^2$ of the variance $\sigma^2$ of the input noise $w(t)$ using the procedure described in \cite[Sec.~V]{bhaskar2013atomic}. Next, we define the quantity
\begin{equation*}
	\beta = \left(L_{\xb} + \log\left(\alpha L_{\xb} \right) + \sqrt{2L \log\left( \alpha L_{\xb} \right)} + \sqrt{\frac{\pi L_{\xb}}{2}} + 1\right)^{\frac{1}{2}}
\end{equation*}
where $\alpha = 8\pi n \log n$. Then, we set
\begin{equation}\label{lambda_value_in_sims}
	\lambda = \hat\sigma \left(1 + \frac{1}{\log n}\right)^{\frac{1}{2}} \beta.
\end{equation}

\subsection{Simulation results}

In this experiment, we generate the signal $y$  of length $L=117$ containing $m=3$ cisoids according to \eqref{signal_model}. The frequencies $\{\theta_k\}_{k=1}^3$ are set as follows. We choose a parameter $\theta_0 \in \{1.1, 1.3, 1.5, \dots, 3.9\}$ and take
$\theta_1=\theta_0-2\Delta_\text{FFT}$, $\theta_2=\theta_0$, and $\theta_3=\theta_0+2\Delta_\text{FFT}$ where $\Delta_\text{FFT}:=2\pi/L$ is the resolution of the FFT method.
In plain words, the separation between two adjacent frequencies is equal to $2\Delta_\text{FFT}$. The complex amplitudes are given by $a_k=e^{i\varphi_k}$ for $k=1, 2, 3$ with phases $\varphi_k$ uniformly distributed in $[0, 2\pi]$. Finally, we add complex Gaussian white noise $w(t)$ to achieve specific SNR levels of $-3, 0, 3$, and $6$ dB, defined as $10\log_{10}(1/\sigma^2)$ using the fact that $|a_k|=1$ for all $k=1, 2, 3$.

\begin{figure}
			\centering
		\begin{subfigure}[b]{.48\columnwidth}
			\includegraphics[width=\linewidth]{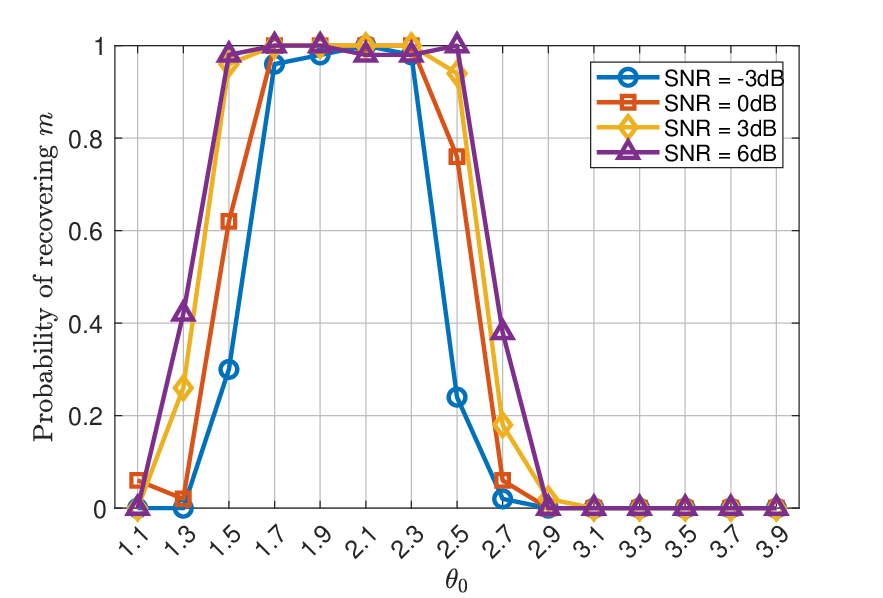}
			\caption{The G-filter $G^{(1)}$}
			\label{subfig:recov_G-filter1}
		\end{subfigure}
		\begin{subfigure}[b]{.48\columnwidth}
			\includegraphics[width=\linewidth]{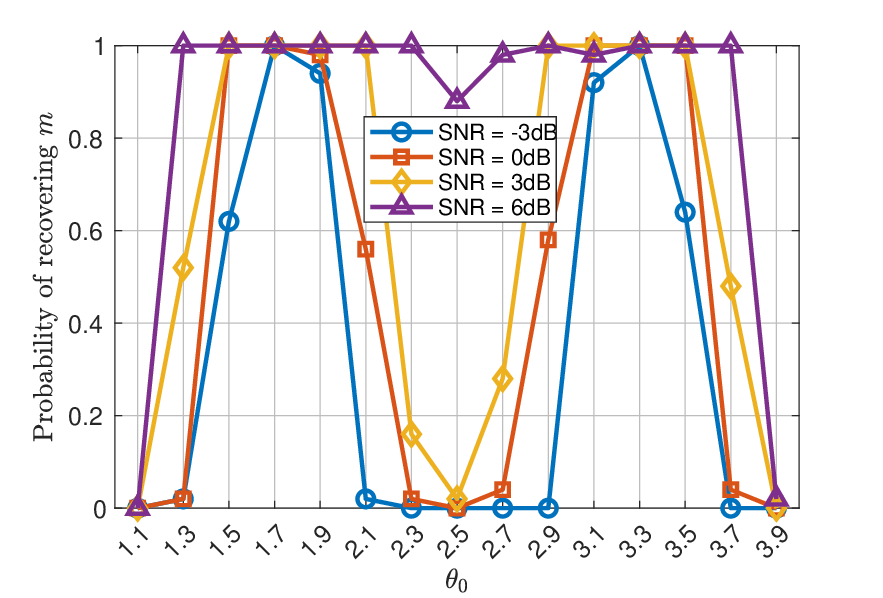}
			\caption{The G-filter $G^{(2)}$}
			\label{subfig:recov_G-filter2}
		\end{subfigure}
        
        \caption{Panels (a) and (b): Probabilities of successfully recovering the number $m=3$ of cisoids versus $\theta_0\in \{1.1, 1.3, \dots, 3.9\}$ with two G-filters $G^{(1)}$ and $G^{(2)}$, respectively. Notice that the markers represent computed values and lines have no meaning.}
		\label{fig_probab_recov_pole}
	\end{figure}

\begin{figure}
		\begin{subfigure}[b]{.48\columnwidth}
			\includegraphics[width=\linewidth]{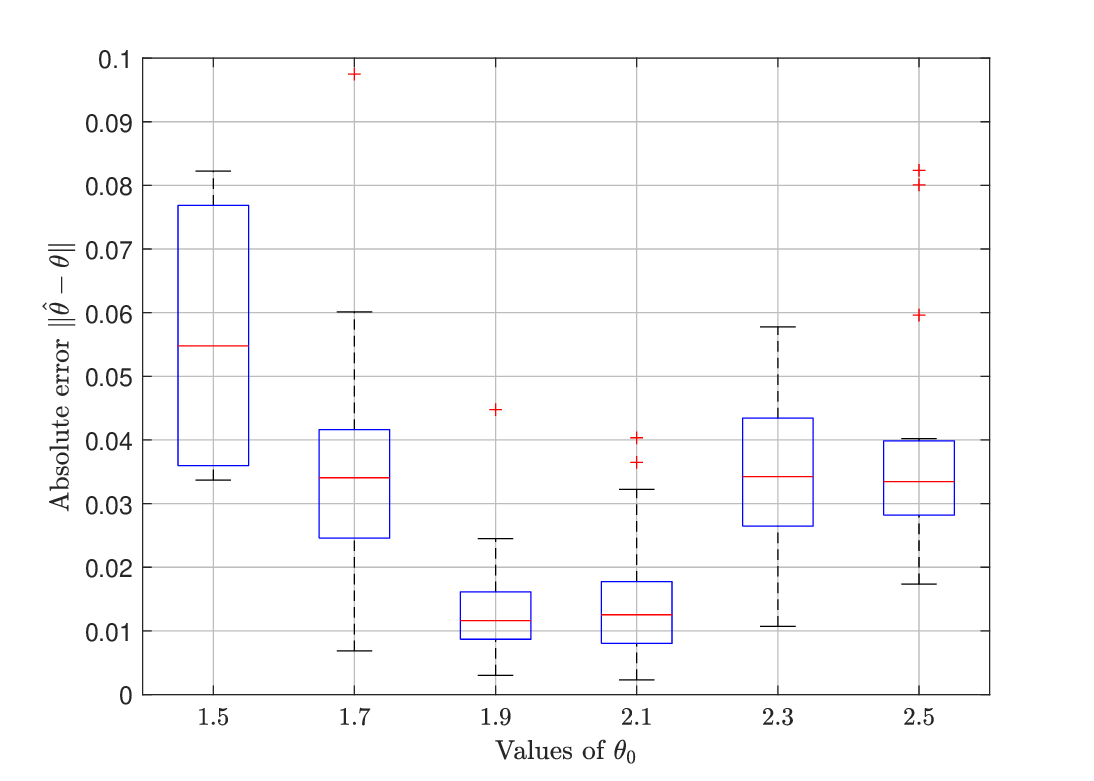}
			\caption{$\SNR=-3$ dB.}
			\label{subfig:poleone_error_SNR1}
		\end{subfigure}
		\hfill
		\begin{subfigure}[b]{.48\columnwidth}
			\includegraphics[width=\linewidth]{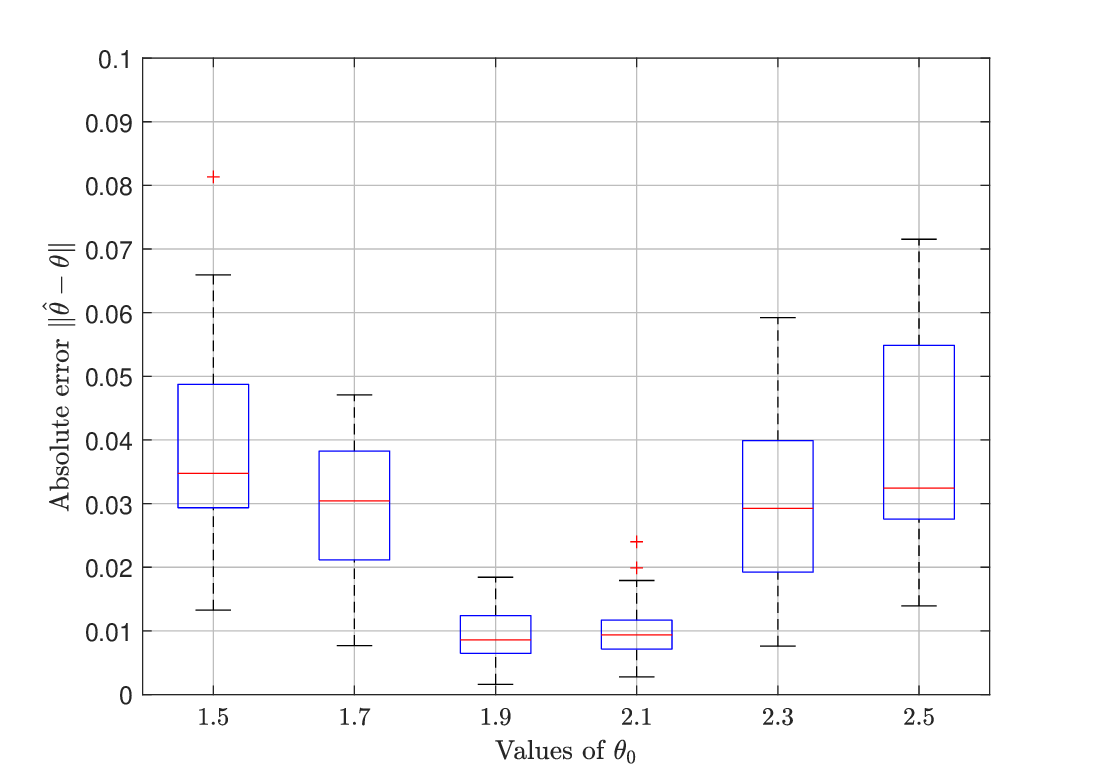}
			\caption{$\SNR=0$ dB.}
			\label{subfig:poleone_error_SNR2}
		\end{subfigure}
		
		\begin{subfigure}[b]{.48\columnwidth}
			\includegraphics[width=\linewidth]{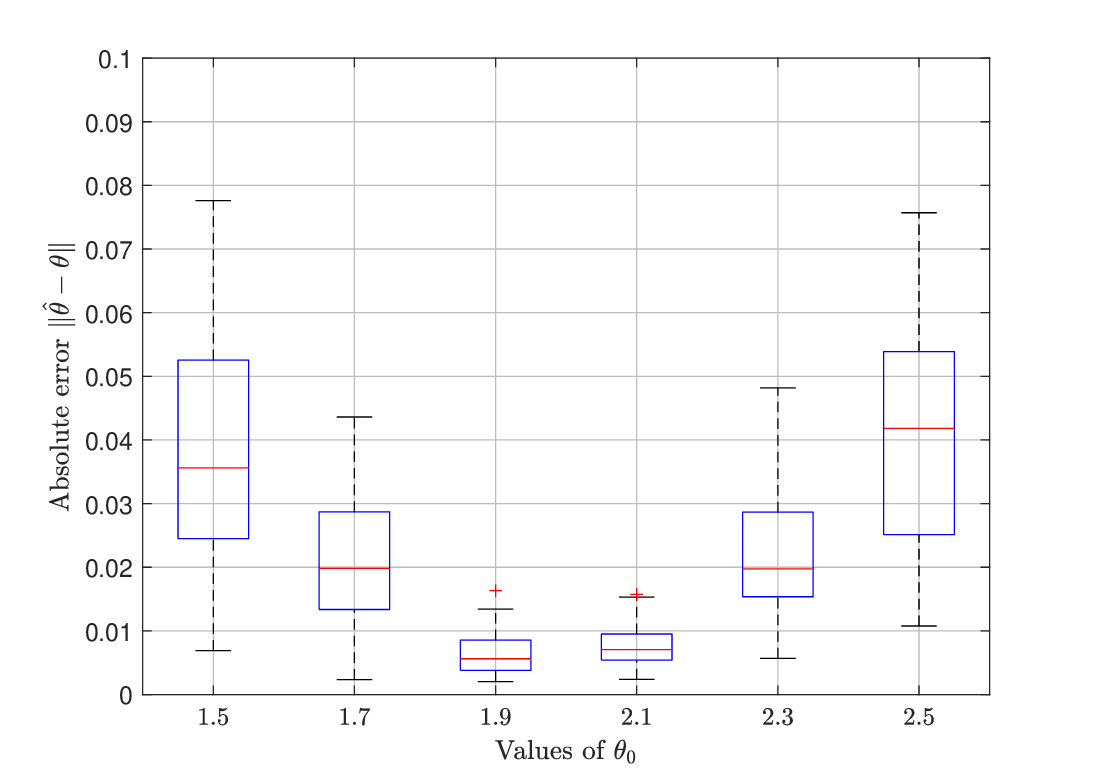}
			\caption{$\SNR=3$ dB.}
			\label{subfig:poleone_error_SNR3}
		\end{subfigure}
		\hfill
		\begin{subfigure}[b]{.48\columnwidth}
			\includegraphics[width=\linewidth]{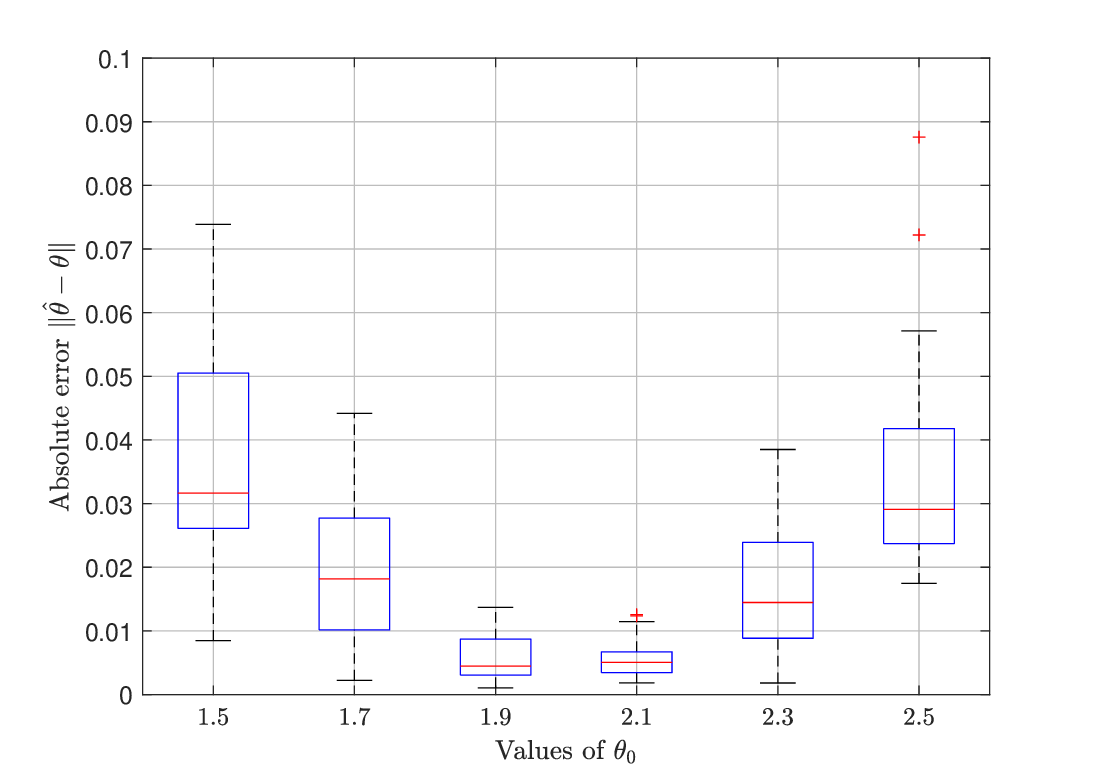}
			\caption{$\SNR=6$ dB.}
			\label{subfig:poleone_error_SNR4}
		\end{subfigure}
		\caption{Panels (a)--(d): Absolute errors $\|\hat{\thetab}-\thetab\|$ of frequency estimation using the G-filter $G^{(1)}$ for $\theta_0^{(1)}\in \{1.5, 1.6, \dots, 2.5\}$ in Subfig.~\ref{subfig:recov_G-filter1} under different SNRs.}
		\label{fig:theta_errors_Gfilter1}
	\end{figure}

\begin{figure}
		\begin{subfigure}[b]{.48\columnwidth}
			\includegraphics[width=\linewidth]{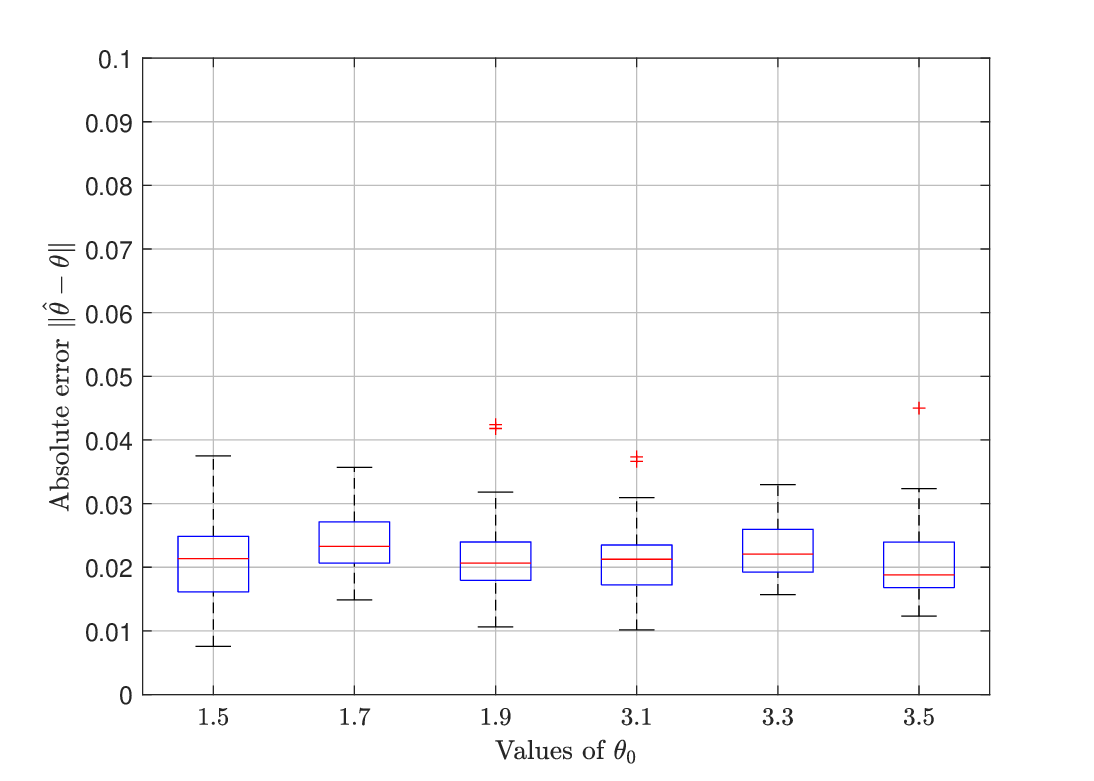}
			\caption{$\SNR=-3$ dB.}
			\label{subfig:poletwo_error_SNR1}
		\end{subfigure}
		\hfill
		\begin{subfigure}[b]{.48\columnwidth}
			\includegraphics[width=\linewidth]{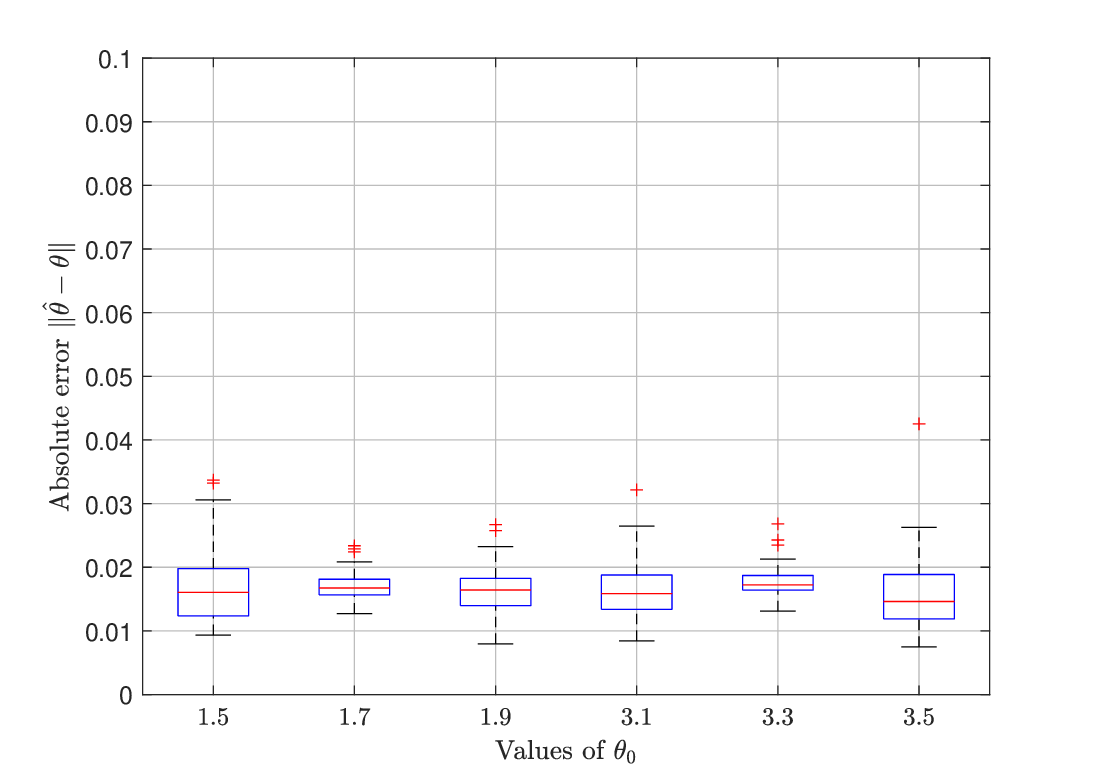}
			\caption{$\SNR=0$ dB.}
			\label{subfig:poletwo_error_SNR2}
		\end{subfigure}
		
		\begin{subfigure}[b]{.48\columnwidth}
			\includegraphics[width=\linewidth]{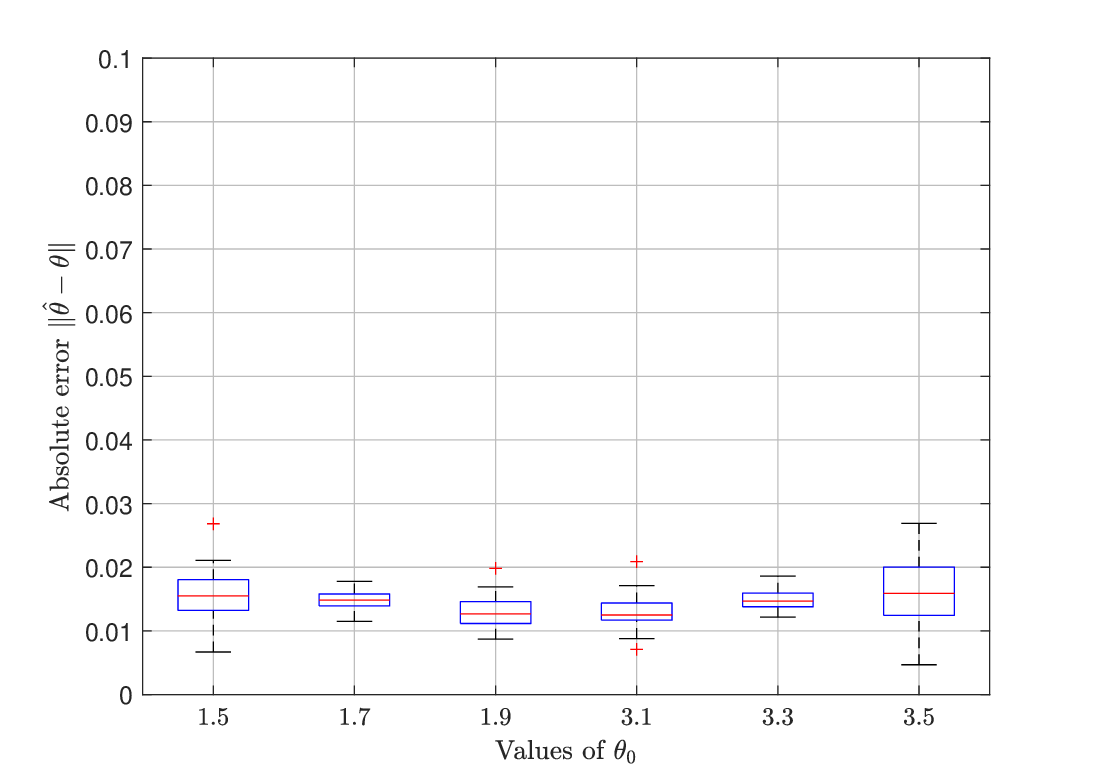}
			\caption{$\SNR=3$ dB.}
			\label{subfig:poletwo_error_SNR3}
		\end{subfigure}
		\hfill
		\begin{subfigure}[b]{.48\columnwidth}
			\includegraphics[width=\linewidth]{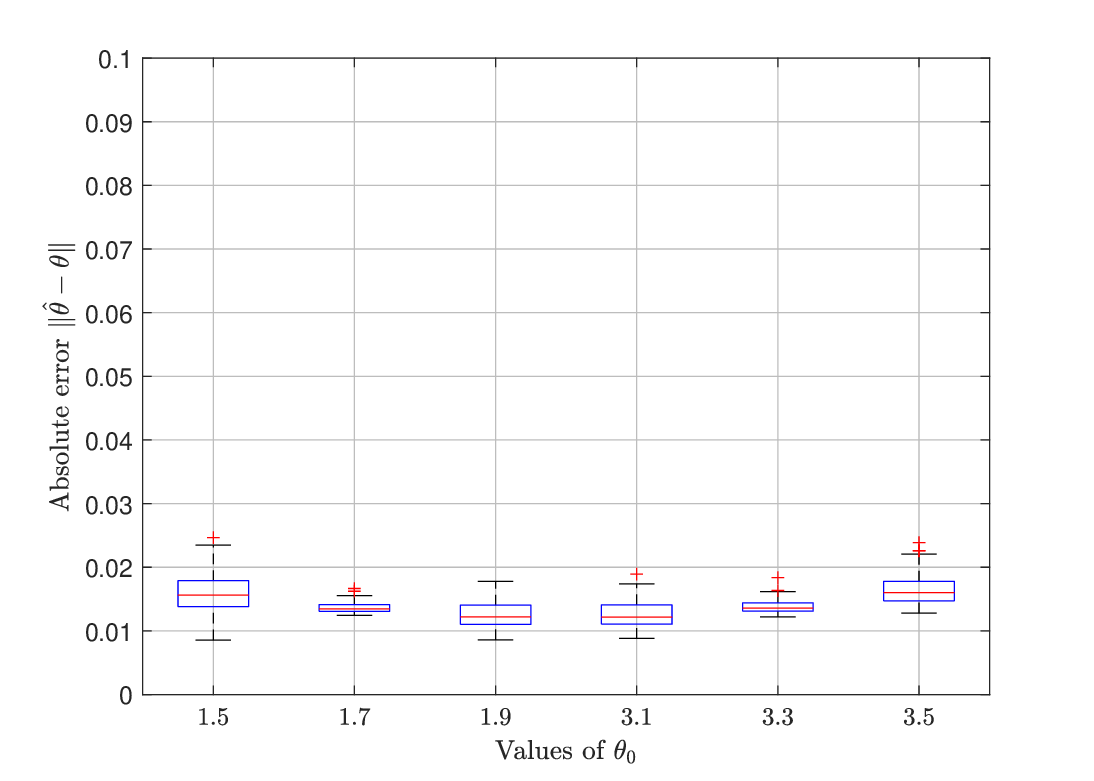}
			\caption{$\SNR=6$ dB.}
			\label{subfig:poletwo_error_SNR4}
		\end{subfigure}

        \caption{Panels (a)--(d): Absolute errors $\|\hat{\thetab}-\thetab\|$ of frequency estimation using the G-filter $G^{(2)}$ for $\theta_0^{(2)}\in \{1.5, 1.7, 1.9, 3.1, 3.3, 3.5\}$ in Subfig.~\ref{subfig:recov_G-filter2} under different SNRs.}
		\label{fig:theta_errors_Gfilter2}
	\end{figure}

    We employ the two G-filters $G^{(1)}$ and $G^{(2)}$ with size $n=20$ in Sec.~\ref{sec_Signal_model} whose squared norms are shown in Figs.~\ref{fig:filter-gain-one} and \ref{fig:filter-gain-two}, respectively, and evaluate the performance of our method by Monte Carlo simulations for those values of $\theta_0$ in the candidate set above. Each Monte Carlo simulation consists of $50$ repeated trials. Following Remark~\ref{rem:input_filtering_revised} for filtering finite-length signals with a parameter $\varepsilon=10^{-3}$, we discard the initial $L_{\srm}^{(1)}=97$ and $L_{\srm}^{(2)}=61$ transient samples, respectively, for $G^{(1)}$ and $G^{(2)}$. This leaves $L_{\xb}^{(1)}=20$ and $L_{\xb}^{(2)}=56$ output vectors, respectively, to construct the data matrix $X$ for the MOV-based SDP \eqref{noisy_SDP}.
    After choosing the regularization parameter $\lambda$ as per Subsec.~\ref{subsec_choice_lambda}, we can solve the SDP  in MATLAB using the CVX \cite{cvx,gb08}. Let $\hat{\tau}_1 \ge \hat{\tau}_2 \ge \dots \ge \hat{\tau}_n$ be the sorted eigenvalues of the optimal $\hat{\Sigma}$. The numerical rank $\hat{r}$ is determined as the smallest index $k$ satisfying $\hat{\tau}_{k+1} < 10^{-3}$ or $\hat{\tau}_k / \hat{\tau}_{k+1} > 10^3$. Then the probability of successfully recovering the true number of $m$ cisoids is computed as
     $P_{\mathrm{succ}} = \frac{1}{50} \sum_{j=1}^{50} \mathbb{I}(\hat{r}_j = m)$, 
where $\hat{r}_j$ denotes the numerical rank of the optimal $\hat{\Sigma}$ in the $j$-th trial of a Monte Carlo simulation and $\mathbb{I}(\cdot)$ represents the indicator function (taking the value $1$ if the condition in the parenthesis holds and $0$ otherwise). The probabilities for the two G-filters $G^{(1)}$ and $G^{(2)}$ when $m=3$ are shown in  Fig.~\ref{fig_probab_recov_pole}, and we see that the recovery probability remains high within the selected passbands even when $\SNR=-3$ dB, demonstrating a robust estimator of the number of cisoids. Next, for each successful trial in the Monte Carlo simulations above, we calculate the absolute error of frequency estimation $\|\hat{\thetab}-\thetab\|$ where $\thetab=(\theta_1, \theta_2, \theta_3)$ denotes the true frequency vector and $\hat\thetab$ the corresponding estimate. Specifically, we evaluate the results for the values of $\theta_0^{(1)} \in \{1.5, 1.6, \dots, 2.5\}$ and $\theta_0^{(2)} \in \{1.5, 1.7, 1.9, 3.1, 3.3, 3.5\}$, which are located inside or near the passbands of the two G-filters $G^{(1)}$ and $G^{(2)}$, respectively. Figs.~\ref{fig:theta_errors_Gfilter1} and \ref{fig:theta_errors_Gfilter2} illustrate the distribution of these errors via the $\mathtt{boxplot}$ in MATLAB. There is a general trend that errors decrease as the SNR increases. In addition, the smallest errors shown in Fig.~\ref{fig:theta_errors_Gfilter1} are observed at $\theta_0=1.7$ and $1.9$ which are close to the center of the selected frequency band $\Ical_1=[1.75, 2.25]$.

\begin{figure}[htbp]
    \begin{subfigure}[b]{.49\columnwidth} 
        \includegraphics[width=\linewidth]{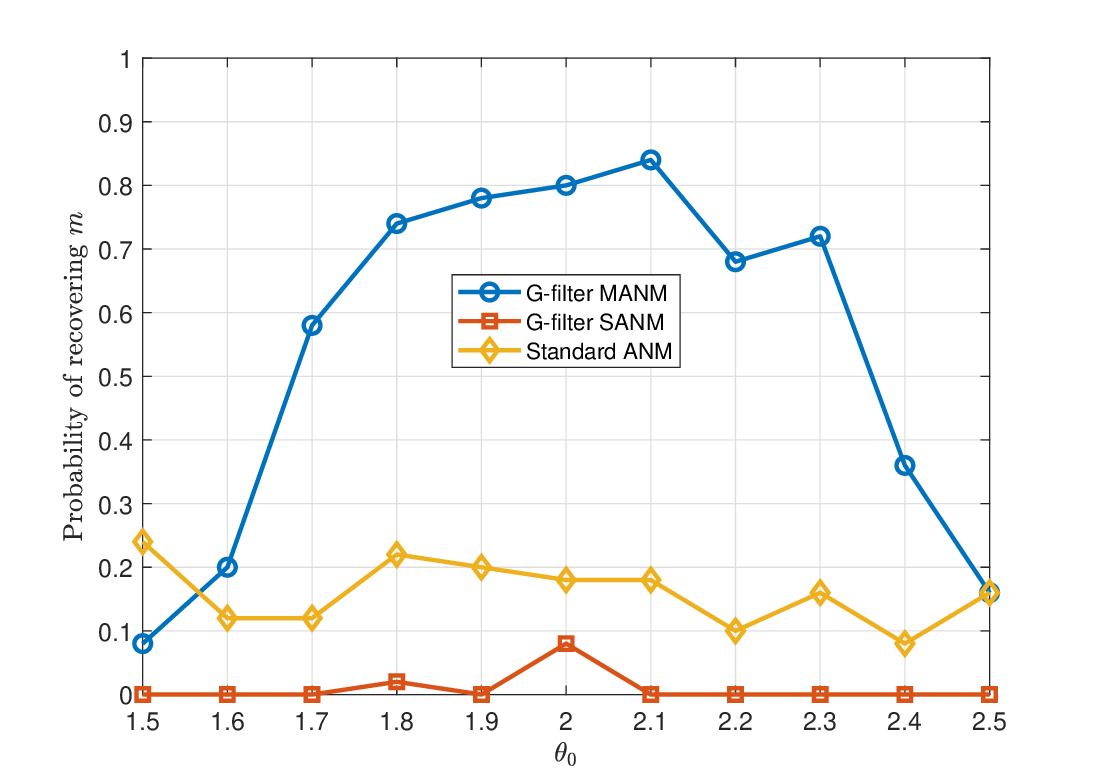} 
        \caption{$\SNR=-3$ dB.}
        \label{subfig:MANM_ANMsnr1}
    \end{subfigure}
    \hfill
    \begin{subfigure}[b]{.49\columnwidth}
        \includegraphics[width=\linewidth]{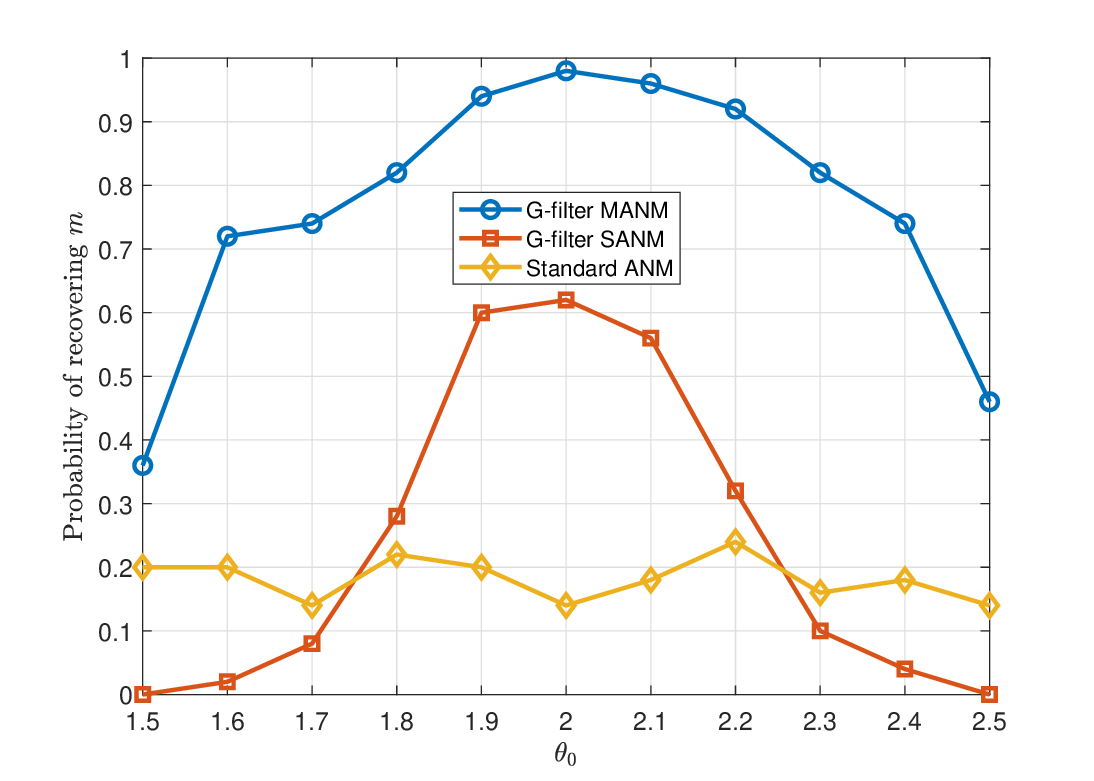}
        \caption{$\SNR=0$ dB.}
        \label{subfig:MANM_ANMsnr2}
    \end{subfigure}
    \hfill   
    \begin{subfigure}[b]{.49\columnwidth}
        \includegraphics[width=\linewidth]{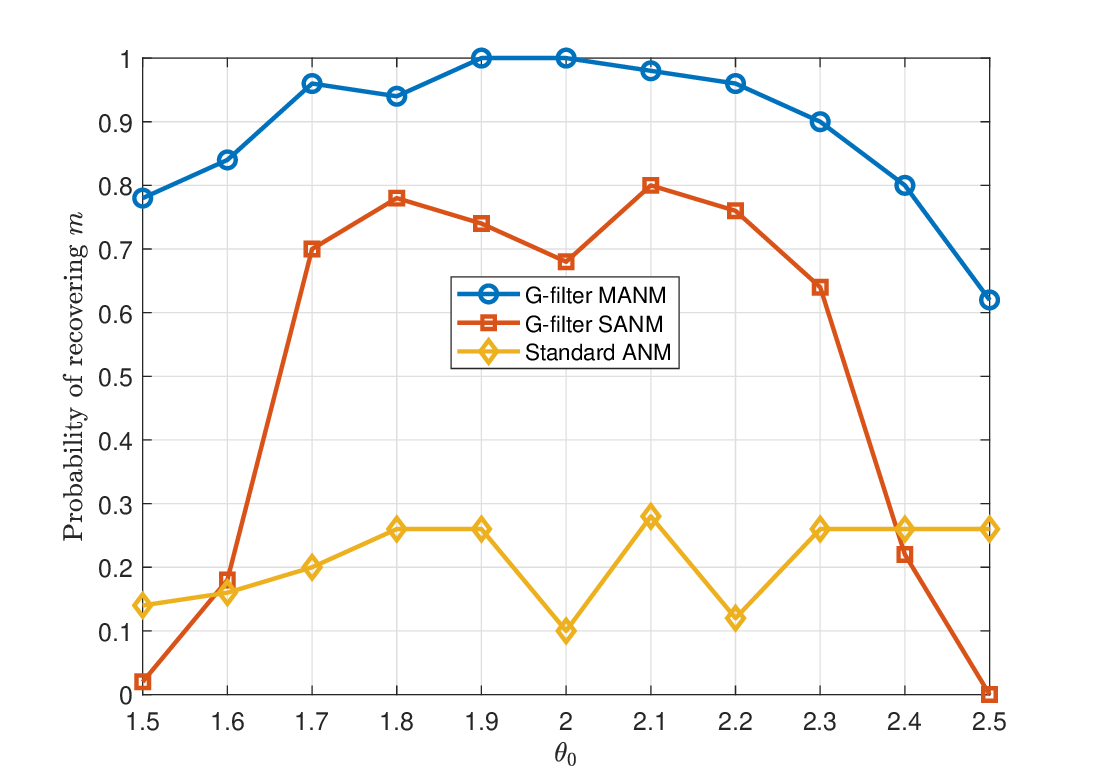}
        \caption{$\SNR=3$ dB.}
        \label{subfig:MANM_ANMsnr3}
    \end{subfigure}
    \hfill 
    \begin{subfigure}[b]{.49\columnwidth}
        \includegraphics[width=\linewidth]{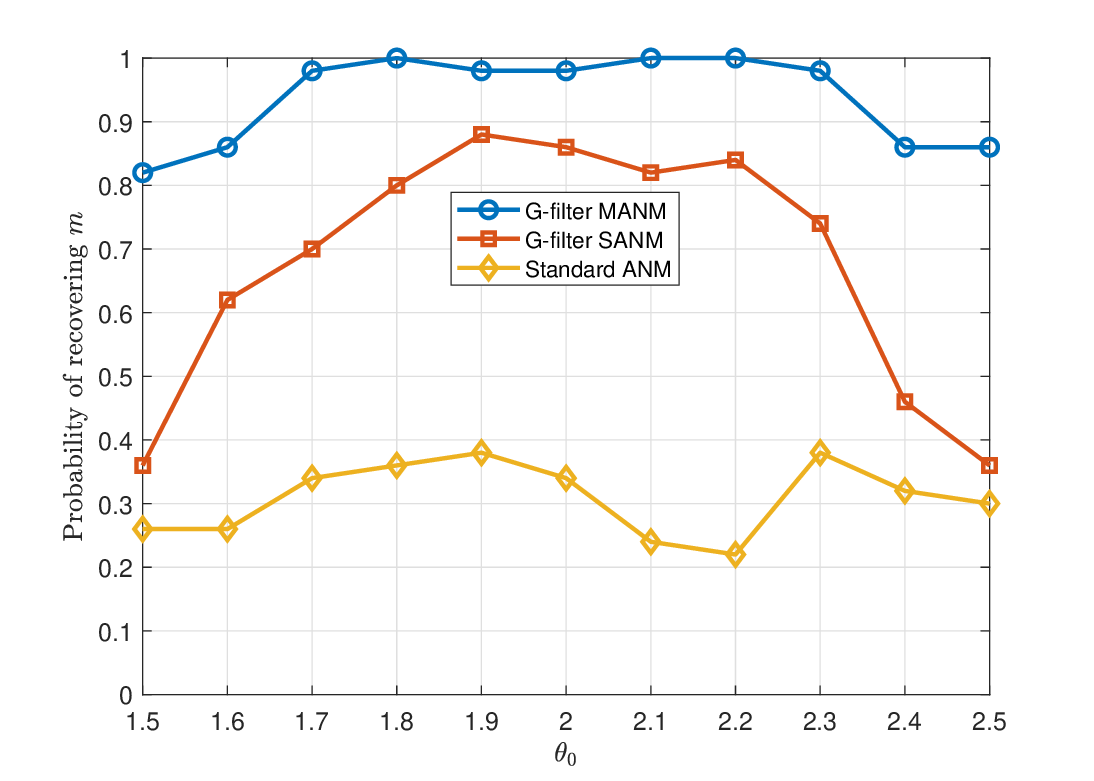}
        \caption{$\SNR=6$ dB.}
        \label{subfig:MANM_ANMsnr4}
    \end{subfigure}
    
    \caption{Panels (a)--(d): Probabilities of successfully recovering the number $m = 2$ of cisoids with G-filter MANM (our approach), the G-filter SANM, and the standard ANM in each Monte Carlo simulation under different SNRs as indicated by captions of the subfigures. Notice that the markers represent computed values and lines have no meaning.}
    \label{fig_MANM_FSANM}
\end{figure}

Furthermore, we evaluate our G-filter MANM approach against the G-filter SANM  \cite{ZHU2026112759} and the standard ANM \cite{bhaskar2013atomic}. Specifically we want to test the \emph{resolution} of the three methods, i.e., their abilities to distinguish two closely located frequencies in the cisoidal signal. To this end, let us consider a signal $y$ in \eqref{signal_model} of length $L=117$ containing $m=2$ cisoids with frequencies $\theta_{1,2} = \theta_0 \pm \Delta_\text{FFT}/2$, where the parameter $\theta_0 \in \{1.5, 1.6, \dots, 2.5\}$. This implies that the separation between the two frequency components $\theta_{1}$ and $\theta_{2}$ is exactly equal to the resolution limit of the FFT method. The complex amplitudes are set as $a_k = 2e^{i\varphi_k}$ with the phase $\varphi_k \sim U[0, 2\pi]$ for $k=1, 2$. The $\SNR$ is calculated as $10\log_{10}(2^2/\sigma^2)$ dB, and we do simulations under different SNRs of $-3, 0, 3$, and $6$ dB. Finally, we adopt the same G-filter $G^{(1)}$ as in Fig.~\ref{fig:filter-gain-one} with size $n=20$, which selects the frequency band $\Ical_1$ and leaves $L_{\xb}=20$ output vectors after truncating the transient samples. Again, we run Monte Carlo simulations of $50$ repeated trials for each choice of $\theta_0$ and the $\SNR$. The probabilities of successfully recovering the number $m=2$ of cisoids for our approach, the G-filter SANM, and the standard ANM are shown in Fig.~\ref{fig_MANM_FSANM}. We observe that our approach outperforms the other two approaches in the recovery probability across the entire frequency range $1.5\leq \theta_0\leq 2.5$ of the simulation under the four values of the SNR (with the only exception of $\theta_0=1.5$ and $\SNR=-3$dB). While resolving closely spaced frequencies proves difficult for the G-filter SANM in the low SNR regime ($-3$ dB), the G-filter MANM approach 
exhibits more robustness against noise, 
demonstrating superior resolution capabilities in the selected band.
    
\section{Conclusion}\label{sec_conc}

In this paper, we propose an improved line spectral estimator that combines Georgiou's filter bank (G-filter) with atomic norm minimization (ANM) using multiple output vectors (MOV) of the filter. By formulating the estimation problem as an SDP and utilizing the Carath\'{e}odory--Fej\'{e}r-type decomposition, our method automatically estimates the number of cisoids (spectral lines) in the signal. Simulations show that exploiting MOV significantly enhances the resolution of the frequency estimate and robustness against noise compared to the standard ANM and the G-filter version of ANM with a single output vector. Future work will focus on developing fast scalable algorithms for (approximate) solutions to the SDPs.

\bibliographystyle{IEEEtran}
\bibliography{ref}

\end{document}